\documentclass[twoside,draft,reqno]{amsart}

\usepackage{latexsym,amscd,amssymb,amsopn,amsthm,amsfonts,amsmath}
\usepackage{color}

\newtheorem{theorem}{Theorem}
\newtheorem{corollary}{Corollary}
\newtheorem{lemma}{Lemma}
\newtheorem{problem}{Problem}
\newtheorem{example}{Example}

\title[On a mixed problem for the parabolic Lam\'e type operator]{On a mixed problem \\ 
 for the parabolic Lam\'e type operator\footnote{This is the preprint version of the paper 
published in J. Inv. Ill-posed Problems, V. 23:6 (2015), 555-570, DOI 10.1515/jiip-2014-0043.}}

\begin{document}

\author{R. Puzyrev}

\address[Roman Puzyrev]
        {Siberian Federal University
\\
         Institute of Mathematics and Computer Science
\\
         pr. Svobodnyi 79
\\
         660041 Krasnoyarsk
\\
         Russia}

\email{effervesce@mail.ru}

\author{A. Shlapunov}

\address[Alexander Shlapunov]
        {Siberian Federal University
\\
         Institute of Mathematics and Computer Science
\\
         pr. Svobodnyi 79
\\
         660041 Krasnoyarsk
\\
         Russia}
         
\email{ashlapunov@sfu-kras.ru}

\begin{abstract}
We consider a boundary value problem for the parabolic Lam\'e type operator being a linearization 
of the Navier-Stokes' equations for compressible flow of Newtonian fluids.
It consists of recovering a vector-function, satisfying the parabolic Lam\'e type system 
in a cylindrical domain, via its values and the values of the boundary stress tensor  
on a given part of the lateral surface of the cylinder. We prove that the problem is ill-posed in 
the natural spaces of smooth functions and in the corresponding  H\"older spaces; besides, 
additional initial data do not turn the problem to a well-posed one. Using the 
Integral Representation's Method we obtain the Uniqueness Theorem and solvability conditions 
for the problem.
\end{abstract}


\maketitle 

\section*{Introduction}

Let, as usual, $\Delta_n$ be the Laplace operator, $\nabla_n$ be the gradient operator and 
$\mbox{div}_n$ be the divergence operator in ${\mathbb R}^n$, $n\geq 2$.  The Navier-Stokes' 
equations for compressible flow of Newtonian fluids over the four-dimensional domain ${\mathcal 
D} \subset {\mathbb R}^3_x \times {\mathbb R}_t $ under the action of a body force $F (x,t)= (F_1 
(x,t),F_2 (x,t),F_3 (x,t))$  can be written in the following form (see \cite[\S 15, 
formulas (15.5), (15.6)]{LaLi}):
\begin{equation} \label{eq.NS}
\rho \Big(\frac{\partial v }{\partial t}  + v \cdot \nabla_3 v\Big) + \nabla_3 p -  \mbox{div}_3 
\Big(\mu_1 \nabla_3 v \Big) - \nabla_3 \Big( \Big(\frac{\mu_1}{3} + \mu_2\Big)  \mbox{div}_3 
v\Big) -  a v = F ,
\end{equation}
where $v (x,t)= (v_1 (x,t), v_2 (x,t), v_3 (x,t))$ is the flow velocity, $\rho (x,t)$ is the 
fluid density, $p (x,t)$ is the pressure, $\mu_j (x,t)$  are (positive) viscosity coefficients 
and  
$$
av =  [(\nabla _3 \mu_1)^* \otimes \nabla_3  - (\nabla_3 \mu_1) \mbox{div}_3] v
$$
is the linear first order summand with $M_1^*$ being the adjoint matrix for a matrix $M_1$ and 
$M^*_1\otimes M_2$ being the Kronecker product of matrices  $M^*_1$ and $M_2$. If the boundary $\partial D$ of $D$ is piece-wise smooth then 
the boundary conditions for this system often involve the force $  \nu p  - \sigma' v  
$ acting on the unit surface area where the force friction (or the boundary viscosity tensor) 
$\sigma'$ has the components
$$
\sigma'_{i,j}   =  \delta_{i,j} \mu_1 \sum_{k=1}^n \nu_k \frac{\partial }{\partial _k}  + 
\mu_1 \nu_j \frac{\partial }{\partial x_i}  + 
(\mu_2 -2\mu_1/3)\nu_i \frac{\partial }{\partial x_j}  , \ 1\leq i,j \leq 3, 
$$
with  $\nu = (\nu_1, \nu_2, \nu_3, \nu_4)$ being the unit normal vector to the surface 
$\partial D$ and $\delta_{i,j}$ being the Kronecker symbol (see \cite[\S 15, formula 
(15.12)]{LaLi}). 

Since the density $\rho$ is positive, a proper linearization of the substantial derivative  term 
$v \cdot \nabla_3 v$  turns (\ref{eq.NS}) into a parabolic Lam\'e type system related to an 
unknown vector $u$: 
$$
L_4 u   = \frac{\partial u }{\partial t} - {\mathcal L} _3  u -  \sum_{j=1}^3 a_j 
(x,t)\frac{\partial u}{\partial x_j} - a_0 (x,t) u  = f
$$
where $a_j (x,t)$, $0\leq j \leq 3$, are $(3\times 3)$ matrices with functional 
entries and 
$$
{\mathcal L}_n  =   \mbox{div}_n \Big(\mu  \nabla_n  \Big) +  
\nabla_n \Big( (\mu + \lambda )  \mbox{div}_n \Big), \, n \geq 2
$$
is the strongly elliptic (with respect to the space variables) formally self-adjoint 
Lam\'e type operator with the Lam\'e coefficients satisfying 
$$
\mu (x,t)>0, \, \, (\mu (x,t)+ \lambda (x,t))\geq 0.
$$
 The regularity \textcolor{red}{(smoothness)} of the Lam\'e  coefficients and the 
matrices $a_j (x,t)$ depends upon the regularity of the density $\rho$ and the viscosity 
coefficients $\mu_j$.

Note that if $\mu$ is constant, $\lambda+\mu =0$ and $a_j=0$, $0\leq j \leq 3$, 
then $L_{4}$ reduces to the heat operator, though, of course, it is known that the heat equation 
is not ideal to model the process of the heat conduction. 

Let $\Omega$ be a bounded domain (i.e. bounded open connected set) in $n$-dimensional  
real space ${\mathbb R}^n$ with the coordinates $x=(x_1, \dots , x_n)$. As usual we denote by   
$\overline{\Omega}$ the closure of $\Omega$, and we denote  by $\partial\Omega$ its boundary.
In the sequel we assume that $\partial \Omega$ is piece-wise smooth. 
As $\partial\Omega$ is piece-wise smooth, the normal 
vector $\nu=(\nu_1, ..., \nu_n)$ is defined almost everywhere on $\partial \Omega$ and 
satisfies $\sum_{j=1}^n \nu_j^2 \ne 0$. 

Let $\Omega_T=\{x\in\Omega,\ 0<t<T\}$  an open cylinder, having the altitude $0<T \leq +\infty$ 
and the base $\Omega$, in $(n+1)$-dimensional  real space ${\mathbb R}^{n+1}= {\mathbb R}^n  
\times \{-\infty<t<+\infty\} $. Let also $\Gamma\subset 
\partial\Omega$ be a non empty connected open (in the topology of $\partial \Omega$) subset 
of $\partial\Omega$ and  $\Gamma_T = \Gamma\times(0,\ T)$. 

In the present paper we consider a mixed boundary problem for the parabolic system 
in the cylindrical domain $\Omega_T$
 \begin{equation*} 
L_{n+1}    = \frac{\partial  }{\partial t} - {\mathcal L} _n -  Au,    
\end{equation*}
where
$$
Au = \sum_{j=1}^n a_j (x,t) 
\frac{\partial  }{\partial x_j} + a_0 (x,t), 
$$
the Lam\'e coefficients and the entries of the $(n\times n)$-matrices 
$a_j (x,t)$, $0\leq j \leq n$, are \textcolor{red}{$C^\infty$-smooth  
in a neighborhood of $\overline{\Omega_T}$ and real analytic with respect to the space 
variables  in a neighborhood of $\overline \Omega$}.

Instead of classical boundary value problems for parabolic equations (see, for instance, 
\cite{lad}, \cite{frid}, \cite{eid}, \cite{sol}) we consider the ill-posed problem, consisting 
in finding a vector-function  satisfying the corresponding parabolic equation in the cylinder 
via its  values  and the values of the boundary stress  tensor with the components
\begin{equation} \label{eq.stress.tensor}
\sigma _{i,j}   =    \mu  \, \delta_ {i,j} \sum_{k=1}^n \nu_k \frac{\partial }{\partial x_k}  + 
\mu \, \nu_j \frac{\partial }{\partial x_i}  + 
\lambda \, \nu_i \frac{\partial }{\partial x_j}  , \ 1\leq i,j \leq n.
\end{equation}
on the given part $\Gamma_T$ of the lateral surface of the cylinder $\Omega_T$ (cf. 
\cite{lor}).  

Using parabolic potentials we prove Uniqueness Theorem and obtain solvability conditions  for the 
problem (cf. \cite{PuSh11} related to  similar results for the heat equation).  Actually, the 
approach was invented  for the investigation of the famous ill-posed Cauchy problem for elliptic 
equations (see, for instance, \cite{AKy} for the Cauchy-Riemann operator,  \cite{Sh96Z} for the 
elliptic Lam\'e operator and \cite{ShTaLMS}, \cite{SheShl}, \cite{Tark36},  for general 
systems with injective principal symbols). 

\section{Preliminaries
}
 
As usual, for $s \in {\mathbb Z}_+$ (here ${\mathbb Z}_+ = {\mathbb N} \cup \{0\}$) and an 
open subset  $D \subset {\mathbb R}^m$ we denote 
$C^s(D)$ the set of all $s$ times continuously differentiable functions in $D$. 
The standard topology of this metrisable space induces uniform convergence on compact 
subsets in $D$ together with all the partial derivatives 
up to order  $s$. 

For $S \subset \partial D$ we denote $C^{s}(D \cup S)$ the set of such functions from the space 
$C^{s}(D)$ that all their derivatives up to order  $s$ can be extended continuously 
onto $D \cup S$. The standard topology of this metrisable space induces uniform convergence on 
compact subsets in $D \cup S$ together with all the partial derivatives 
up to order  $s$. In particular, for bounded domains,  $C^{s}(D \cup \partial 
D)= C^{s}(\overline D)$ is a Banach space. 
If $D$ is an unbounded  then the Banach 
space $C^s_b (D\cup S)$ consists of $s$-times differentiable functions in $D \cup S$ 
with bounded derivatives up to order $s$ and it is endowed with the standard $sup$-norm.  
Then $C^s_b (\overline D) = C^s (\overline D)$ for a bounded domain $D$. 

Apart from the standard functional spaces, we need also spaces taking into account the specific 
properties of parabolic equations in ${\mathbb R}^{n+1}= {\mathbb R}^n \times  
\{-\infty<t<+\infty\}$. Namely, let $C^{1,0} (\Omega_T)$ be the set of continuous 
functions $u$ in  $\Omega_T$, having in $\Omega_T$ continuous partial derivatives 
$u_{x_i}$, and let $C^{2,1} (\Omega_T)$ denote  the set of continuous 
functions in  $\Omega_T$, having in $\Omega_T$ continuous partial derivatives 
$u_{x_i}$, $u_{x_ix_j}$, $u_t$. The standard topology of this metrisable space induces 
uniform convergence on compact 
subsets in $D$ together with all the partial derivatives used in its definition. 

As before, for $S \subset \partial \Omega_T$ we denote by $C^{1,0}(\Omega_T \cup S)$ the set of 
such functions $u$ from the space $C^{1,0}(\Omega_T)$ that their derivatives  $u_{x_i}$ 
can be extended continuously onto $\Omega_T \cup S$. The standard topology of this 
metrisable space induces uniform convergence on compact subsets of $\Omega_T \cup S$ of both the 
functional sequences and the corresponding sequences of first partial derivatives ${x_i}$.
Clearly,  $C^{1,0}(\Omega_T \cup \partial \Omega_T)= C^{1,0}(\overline {\Omega_T})$ 
is a Banach space. 
Similarly to the standard spaces, if $D \subset {\mathbb R}^{n+1}$
the Banach space $C^{1,0}_b (D \cup S)$ consisting of bounded  $C^{0,1} (D \cup S)$-functions 
with bounded derivatives  $u_{x_i}$ in $D \cup S$ and it is endowed with the norm  
$$
\|u\|_{C^{1,0}_b (D \cup S)}= \sup_{(x,t) \in D \cup S} 
|u(x,t)| + \sup_{(x,t) \in D \cup S} \sum_{j=1}^n
|\frac{\partial u(x,t)}{\partial x_j }|.
$$ 
Then $C^{0,1}_b (\overline D ) = C^{0,1} (\overline D)$ for a bounded domain $D$. 

The space of $n$-vector-functions $u=(u_1, \dots , u_n)$ of a 
class $\mathfrak C $ will be denoted by $[\mathfrak C]^n$.

Let now $\theta$ be such positive constant that for all $(x,t) \in \overline \Omega_T$ we have
$$
\mu(x,t) \geq \theta, \, (\lambda(x,t) + 2\mu(x,t)) \geq \theta. 
$$
Then a direct calculation shows that for all $\zeta \in {\mathbb R}^n$ we have
$$
\det{\Big(\mu (x,t)|\zeta|^2 (\sqrt{-1})^2 I_n+ (\lambda(x,t) + \mu(x,t)) \zeta \zeta^T 
(\sqrt{-1})^2 - \kappa I_n \Big)} = 
$$
$$
(-1)^n (\mu(x,t)|\zeta|^2 + \kappa)^{n-1} ((2\mu (x,t)+\lambda(x,t))|\zeta|^2+\kappa) 
$$
where $I_n$ is the unit $(n\times n)$ matrix and $\zeta^T$ is the transposed vector for $\zeta$.
Hence the roots of this polynomial (with respect to $\kappa$) are 
$$
\kappa_1 (x,t,\zeta)= - (2\mu (x,t)+\lambda(x,t))|\zeta|^2 , \, \kappa_2 (x,t,\zeta) = -\mu 
(x,t)|\zeta|^2 
$$
and, for all $(x,t)\in \overline \Omega_T$, we have  
$$
\max \Big( \sup_{|\zeta|=1} \kappa_1 (x,t,\zeta), \sup_{|\zeta|=1}  \kappa_2 (x,t,\zeta) \Big ) 
\leq - \theta,
$$
i.e. the operator $L_n$ is uniformly parabolic (according to Petrovskii) on $\overline \Omega_T$.

Now we assume that  there is a $n$-dimensional domain $U \supset \overline \Omega$ such 
that the Lam\'e coefficients $\mu (x,t)$, $\lambda (x,t)$ 
and the entries of the $(n\times n)$-matrices $a_j (x,t)$, $0\leq j \leq n$, are 
\textcolor{red}{$C^\infty$-smooth  in $\overline {U_T}$ and real  
analytic with respect to the space variables  in $U$.} 

Under the assumptions, the following properties hold true for parabolic 
operator $L_{n+1}$, which will be crucial for the approach below.

\begin{theorem} \label{t.prop.1}
Each  weak solution $u$ to $L_{n+1} u =0$ in the domain $\Omega_T \subset U_T$ belongs 
to $C^\infty (\Omega_T)$ and it is actually real analytic with respect to variables $x$ in 
$\Omega$. 
\end{theorem}

\begin{theorem} \label{t.prop.2}
The operator $L_{n+1}$ has a fundamental solution in $U_T$, i.e. a 
$(n \times n)$-matrix $\Phi (x,t,y,\tau)$ satisfying 
\begin{equation} \label{eq.fund.1}
(L_{n+1})_{x,t} \Phi (x,t,y,\tau) = 0, \ 
(L^*_{n+1})_{y,\tau} \Phi (x,t,y,\tau) =0, \ \mbox{ if } (x,t) \ne 
(y,\tau), 
\end{equation}
$$
jj
$$
with the formal adjoint operators 
$$ 
(L^*_{n+1})_{y,\tau} = - \frac{\partial  }{\partial \tau}  - ({\mathcal L} _n)_y  - A^* , 
\, A^* = - \sum_{k=1}^n \frac{\partial }{\partial y_k} (a^*_k (y,\tau) \cdot )  + a^*_0 (y,\tau).
$$
\end{theorem}

{\bf Proof} See, for instance, \cite[Ch. 2]{eid}.
\hfill $\square$

We need a sort of an integral representation,  similar to the famous Green Formula for 
the Laplace Operator, constructed with the use the fundamental solutions. 
More precisely, consider the cylinder type domain  
$\Omega_{T_1,T_2} = \Omega_{T_2} \setminus \overline{\Omega_{T_1}}$ and a closed 
measurable set $S \subset \partial \Omega$. 

Let $\sigma$ be the tensor with the components given by (\ref{eq.stress.tensor}) and 
$$
\tilde \sigma = \sigma - \sum_{k=1}^n a^*_k (x,t) \nu_k (x) \textcolor{red}{+}   
[(\nabla_n \mu (x,t))  \nu^* (x)  - \nu (x) (\nabla_n \mu (x,t))^*]    .
$$

For functions $f \in C (\overline{\Omega_{T_1,T_2}})$, $v \in C(S_T)$, $w \in C(S_T)$, $h \in 
C(\overline{\Omega})$ we set  
\begin{equation} \label{3.4}
I_{\Omega, T_1} h (x,t)= \int\limits_{\Omega}\Phi^*(x,t,y,T_1)h(y) dy,  
\end{equation}
\begin{equation} \label{3.1}
G_{\Omega, T_1} f (x,t)=\int\limits_{T_1}^t\int\limits_\Omega \Phi^*(x,t,y,\tau)f(y, \tau)dy 
d\tau, 
\end{equation}
\begin{equation} \label{3.2}
V_{S,T_1} v (x,t)=\int\limits_{T_1}^t\int\limits_{S} \Phi^*(x,t,y, \tau) v(y, \tau)ds(y)d\tau, 
\end{equation}
\begin{equation} \label{3.3}
W_{S,T_1} w (x,t)= - \int\limits_{T_1}^t\int\limits_{S} [\tilde \sigma _y 
\Phi(x,t,y,\tau)]^* w(y, \tau)ds(y) d\tau, 
\end{equation}
where $ds$ is the volume form on $S$ induced from ${\mathbb R}^{n}$. All these functions are 
called \emph{Parabolic Potentials} with densities  $f$, $v$, $w$ and $h$ respectively. In our 
situation these are convergent improper integrals depending on vector parameter $(x,t)$ in 
the neighborhood $U$ of the cylinder $\overline{\Omega_{T_1, T_2}}$ in ${\mathbb R}^{n+1}$  
(see, for instance, \cite[Ch. 4, \S 1]{lad}, \cite[Ch. 3, 
\S 10]{landis}, \cite[Ch. 1, \S 3 and Ch. 5, \S 2]{frid}). 
The potential $I_{\Omega,T_1} (h)$ is sometimes called \emph{Poisson type integral} for the Lam'e 
type Operator, the functions $G_{\Omega,T_1} (f)$,  $V_{S,T_1} (v)$, $W_{S,T_1} (w)$ are often 
referred to as  \emph{Parabolic Volume Potential}, \emph{Parabolic Single Layer Potential} and 
\emph{Parabolic Double Layer Potential} respectively. 


\begin{lemma} \label{l.Green} 
For all $0 \leq  T_1 < T_2 \leq T$ and all  $u \in C^{2,1} (\Omega_{T_1,T_2}) 
\cap C^{1,0} (\overline{\Omega_{T_1,T_2}})$ with 
$L_{n+1}u \in C (\overline{\Omega_{T_1,T_2}})$ the following 
formula holds:
\begin{equation} \label{eq.Green}
 \left( 
I_{\Omega,T_1} u   + G_{\Omega,T_1} L_{n+1}u  + V _{\partial \Omega, T_1} 
\sigma u  +  W_{\partial \Omega, T_1} u \right) (x,t)  = 
\left\{
\begin{aligned}
u(x, t),\ (x,t)\in \Omega_{T_1,T_2}  \\
0,\ (x, t)\in U_T \setminus \overline{\Omega_{T_1,T_2}} .
\end{aligned}
\right.  
\end{equation}
\end{lemma}

{\bf  Proof.}  Indeed, it follows from Gau\ss{}-Ostrogradskii Formula that 
\begin{equation} \label{eq.Green.0}
\int_{\partial \Omega} v^* \sigma u
= \int_{\partial \Omega} v^* ( {\mathcal L}_n u+ au) dy + {\mathfrak D}_\Omega (u,v)
\end{equation}
for all $u,v \in C^{1,0} (\overline \Omega_{T_1,T_2})$ with 
$L_{n+1} u \in C (\overline \Omega _{T_1,T_2})$, 
where 
\begin{equation} \label{eq.Green.a}
au = [(\nabla _n \mu)^* \otimes \nabla_n  - (\nabla_n \mu) \mbox{div}_n] u, 
\end{equation}
\begin{equation} \label{eq.Green.d}
{\mathfrak D}_\Omega (u,v) = \int_\Omega \Big( \mu (\nabla_n v)^*\nabla_n u + \mu 
(\nabla_n v)^*\otimes \nabla_n u) + \lambda (\mbox{div}_n v)^* \mbox{div}_n u  \Big) dy.
\end{equation}
On the other hand, by Gau\ss{}-Ostrogradskii Formula,  
$$
\int_{\partial \Omega}[(\nabla_n \mu (x,t))  \nu^* (x)  - \nu (x) (\nabla_n \mu (x,t))^*] ds (y) 
= \int_\Omega v^* (au - (a^*v)^* u ) \, dy. 
$$
Therefore 
$$
\int_{\partial \Omega} \Big(v^* \sigma u - (\tilde \sigma v )^*u\Big) ds(y) = 
\int_\Omega \Big(v^* ({\mathcal L}_n u + Au) - ({\mathcal L}_nv + A^*v)^* u \Big) dy
$$
for all $u,v \in C^1 (\overline \Omega)$ with ${\mathcal L}_nu, {\mathcal L}_n v \in C (\overline 
\Omega)$. Hence, again by Gau\ss{}-Ostrogradskii Formula, we obtain the (first) Green formula 
for the Lam\'e type operator: 
\begin{equation} \label{eq.Green.1}
\int_\Omega [v^* (y,T_1) u (y,T_1) -  v^* (y,T_2) u (y,T_2)]dy \textcolor{red}{-} \int_{T_1}
^{T_2}\int_{\partial \Omega } \Big(v^* \sigma u - (\tilde \sigma v )^*u\Big)  ds(y) \, dt  = 
\end{equation}
$$
 \int_{\Omega_{T_1,T_2}} \Big(v^* L^*_{n+1} u - (L_{n+1} v)^* u \Big)  dt \,  dy 
$$
for all $u,v \in C^{1,0} (\overline \Omega_{T_1,T_2})$ with 
$L_{n+1} u, L_{n+1} v \in C (\overline \Omega _{T_1,T_2})$.

It follows from the definition of the fundamental solution, that 
$$
(L_{n+1})_{x,t} \Phi (x,t,y,\tau) = \delta (x-y, t-\tau), \ 
(L^*_{n+1})_{y,\tau} \Phi (x,t,y,\tau) =\delta (x-y, t-\tau), 
$$
$$
\Phi (x,t,y,\tau) =0 \mbox{ for } \tau>t.
$$
Then, using the standard argument (see, for instance, \cite[Ch. 6, \S 12]{svesh} 
for the heat equation), we see that Green's Formula (\ref{eq.Green}) follows 
from (\ref{eq.Green.1}) and Fubini Theorem.
\hfill $\square$


\begin{theorem}[Uniqueness Theorem] \label{t.1}
 If $\Gamma$ has at least one interior point (on $\partial\Omega$), 
 and function $u\in C^{2,1}(\Omega_T)\cap C^{1,0}(\Omega_T\cup\overline{\Gamma_T})$ satisfies 
 $L_{n+1} u \equiv 0$ in $\Omega$, $u\equiv 0$ on $\Gamma_T$, $\sigma  u \equiv 0$ on $\Gamma_T$, 
 then $u \equiv 0$ in  $\Omega_T$.  
\end{theorem}

\textbf{Proof.} Under the hypothesis of the theorem there is an interior point $x_0$  
on $\Gamma$. Then there is such a number $r>0$ that $B(x_0,\ r)\cap\partial\Omega\subset\Gamma$ 
where  $B(x_0,\ r)$ is ball in $U \subset {\mathbb R}^n$ with center at $x_0$ and radius  
$r$. Fix an arbitrary point $(x',t') \in \Omega_T$. It is clear that there is a domain $\Omega' 
\ni x'$ satisfying $\Omega ' \subset \Omega$ and $\Omega ' \cap \partial \Omega \subset \Gamma 
\cap B(x_0,\ r)$. Then $(x', t') \in {\Omega'_{T_1,T_2}}$ with some $0 < T_1 < T_2 <T$.

But $u\in C^{2,1}(\Omega'_{T_1,T_2})\cap C^{1,0}(\overline{\Omega'_{T_1,T_2}})$ and  
$L_{n+1}u = 0$ in $\Omega'_{T_1,T_2}$ under the hypothesis of the theorem. Hence 
formula (\ref{eq.Green}) implies:  
\begin{equation} \label{eq.Green.prime}
I_{\Omega',T_1} u (x,t) + V _{\partial \Omega' \setminus \Gamma, T_1} \sigma u  
(x,t)+ W_{\partial \Omega' \setminus \Gamma, T_1} u (x,t)  = 
\left\{ 
\begin{aligned}
u(x, t),\ (x,t)\in \Omega'_{T_1,T_2} , \\
0,\ (x, t)\in U_T \setminus \overline{\Omega'_{T_1,T_2}} , 
\end{aligned}
\right.  
\end{equation}
because  $u\equiv \sigma u \equiv 0$ on $\Gamma_T$.

Taking into account \textcolor{red}{the character of the singularity of the kernel (see 
\cite[Theorem 2.2]{eid})}  $\Phi (x,y,t,\tau)$ we conclude that the following properties are 
fulfilled for the integrals, depending on parameter, from the right hand side of identity 
(\ref{eq.Green.prime}):
\begin{equation*}
 I_{\Omega',T_1} (u)  \in C^{2,1} (U_{T_1,T_2}), 
\end{equation*}
\begin{equation*}
W_{\partial \Omega' \setminus \Gamma, T_1} u, \, V _{\partial \Omega' \setminus \Gamma, T_1} 
\sigma u  \in C^{2,1} ((U \setminus (\partial \Omega' \setminus \Gamma))_{T_1,T_2}) 
\end{equation*}
(see, for instance, \cite[Ch. 4, \S 1]{lad}, \cite[Ch. 3, \S 10]{landis} 
or \cite[Ch. 1, \S 3 and  Ch. 5, \S 2]{frid}).
Moreover, as $\Phi$ is a fundamental solution to Lam\'e  type operator 
then  using (\ref{eq.fund.1}) and Leibniz rule for differentiation of integrals 
depending on parameter we obtain: 
\begin{equation*}
 L_{n+1} I_{\Omega',T_1} u = 0  \mbox{ in }  U_{T_1,T_2}, 
\end{equation*}
\begin{equation*}
L_{n+1}V _{\partial \Omega' \setminus \Gamma, T_1} \sigma u 
= L_{n+1} W_{\partial \Omega' \setminus \Gamma, T_1} u = 0 \mbox{ in  }
(U \setminus (\partial \Omega' \setminus \Gamma))_{T_1,T_2}.
\end{equation*}

Hence the function  
\begin{equation*}
P(x,t) = I_{\Omega',T_1} u (x,t) + 
 V _{\partial \Omega' \setminus \Gamma, T_1} \sigma u  
(x,t)+ W_{\partial \Omega' \setminus \Gamma, T_1} u (x,t)   , 
\end{equation*}
satisfies the Lam\'e type equation 
\begin{equation*}
 (L_{n+1} P) (x,t) = 0  \mbox{ in } (U \setminus (\partial \Omega' \setminus \Gamma))_{T_1,T_2} .
\end{equation*}
This implies that the function $P (x,t)$ is real analytic with respect to the space variable 
$x \in U \setminus (\partial \Omega' \setminus \Gamma)$ for any $T_1<t<T_2$ (see, 
for instance, \cite[Ch. VI, \S 1, Theorem 1]{mih}). 
In particular, by the construction the function $P (x,t)$ is real analytic with respect to $x$ 
in the ball $B(x_0,r)$ and it equals to zero for $x \in B(x_0,R) \setminus \overline \Omega$ for 
all $T_1<t<T_2$. Therefore, the Uniqueness Theorem for real analytic functions yields 
$P(x,t) \equiv 0$ in $(U \setminus (\partial \Omega' \setminus \Gamma))_{T_1,T_2}$, 
and in the cylinder $\Omega'_{T_1,T_2}$, containing the point $(x',t')$. 
Now it follows from (\ref{eq.Green.prime}) that  $u (x',t') =P(x',t') = 0$ and then, since 
the point $(x',t') \in \Omega_T$ is arbitrary we conclude that $u \equiv 0$ in $\Omega_T$.
The proof is complete. 
\hfill $\square$

\begin{example} \label{ex.2}
Let $\mu=1$, $\lambda=-1$ and $a_j=0$, $0\leq j \leq n$. 
Then $L_{n+1}$ reduces to the heat operator: 
$$
L_{n+1} = \frac{\partial}{\partial t}- \Delta_n 
$$ 
and corresponding fundamental solution is given by $\Phi (x,y,t,\tau)= \varphi_0 
(x-y,t-\tau) I_n$  
where 
\begin{equation*}
\varphi_0(x,t)=\begin{cases}
\frac{1}{\left(2\sqrt{\pi t}\right)^n}\ e^{-\frac{|x|^2}{4 t}} & \mbox{ if } t>0,\\
0 & \mbox{ if } t\leq 0.
\end{cases} 
\end{equation*}
In this case $\tilde \sigma = \sigma =\frac{\partial}{\partial \nu} $.
\end{example}

\begin{example} \label{ex.3}
Let $\mu$, $\lambda$ be constant and $a_j=0$, $0\leq j \leq n$. 
Then $L_{n+1}$ reduces to the parabolic Lam\'e operator
$$
L_{n+1} = \frac{\partial}{\partial t}- {\mathcal L}_n 
$$ 
and corresponding fundamental solution $\Phi (x,y,t,\tau)$ 
is given by $(n\times n)$-matrix with components 
$ \Phi_{i,j} (x,y,t,\tau) = \varphi_{i,j} (x-y,t-\tau)$ where 
\begin{equation*}
\varphi_{i,j} (x,t)= 
\varphi_0 (x, \mu t) \delta _{i,j} 
+ \int_{\mu t}^{(2\mu+\lambda) t} \frac{\partial^2 \varphi_0(x,s)}{\partial x_j 
\partial x_i} ds, 
\end{equation*}
(see, for instance, \cite{eid}). In this case $\tilde  \sigma = \sigma = \mu\frac{\partial}
{\partial \nu} +  \mu \, \nu^* \otimes  \nabla_n + \lambda \, \nu \, \mbox{div}_n .$
\end{example}

\section{The boundary problem}

Green formula (\ref{eq.Green}) and the Uniqueness Theorem \ref{t.1} suggest 
us to consider two kind of problems for the parabolic Lam\'e type equation.  

Let vector-functions 
$$u^{(0)}(x)  \in [C(\overline \Omega)]^n, \, 
f(x,t) \in [C_b(\overline{\Omega_T})]^n, 
$$
$$
u^{(1)}(x,t) \in [C^{1, 0} (\overline{\Gamma_T}) \cap C_b(\overline{\Gamma_T})]^n, 
\, u^{(2)}(x,t)  \in [C_b(\overline{\Gamma_T})]^n \, 
$$ 
be given.

\begin{problem} \label{pr.1}
Find 
a vector-function 
$
u(x,\ t) \in [C^{2, 1}(\Omega_T) \cap C^{1,0}(\Omega_T \cup \overline{\Gamma_T})  ]^n 
$ 
satisfying the Lam\'e type  equation 
\begin{equation} \label{1.2}
L_{n+1} u = f \mbox{ in } \Omega_T 
\end{equation}
and boundary conditions 
\begin{equation} \label{2.1}
 u(x,t)=u^{(1)}(x,\ t) \mbox{ on } \overline{\Gamma_T}, 
\end{equation}
\begin{equation} \label{2.2}
\sigma u(x,t) =u^{(2)}(x,\ t) \mbox{ on } \overline{\Gamma_T}. 
\end{equation}
\end{problem}

Note that, if the surface $\Gamma$ and the data of the problem are real analytic then 
the Cauchy-Kovalevsky Theorem implies that  Problem \ref{pr.1}  can not have more than one 
solution in the class of (even formal) power series. However the theorem does not imply 
the existence of solutions to Problem \ref{pr.1} because it grants the solution in a 
small neighborhood of the surface  $\Gamma_T$ only (but not in a given domain $\Omega_T$!).  
In any case, we do not assume the real analyticity of $\Gamma$ and the data $u^{(1)}$, 
$u^{(2)}$  and $f$.

\begin{corollary} \label{c.1}
If $\Gamma$ has at least one interior point (on $\partial\Omega$) then 
Problem \ref{pr.1} has no more that one solution.
\end{corollary}

\textbf{Proof.}
Let  $v(x,\ t)$ and $w(x,\ t)$ be two solutions to Problem \ref{pr.1}. Then 
function $u=(v-w) \in C^{2,1}(\Omega_T)\cap C^{1,0}(\Omega_T\cup \overline{\Gamma_T})\cap 
C(\overline{\Omega_T}\setminus (\partial \Omega \setminus \Gamma)_T)$ is a solution 
to the corresponding problem with $f=0,\ u_1=0,\ u_2=0$. Using  \ref{t.1} we conclude that 
$u$ is identically zero in  $\Omega_T$. 
\hfill $\square$

Thus, the Uniqueness Theorem \ref{t.1} implies that the data of Problem 
\ref{pr.1}  are suitable in order to uniquely define its solution. 

Easily, Problem \ref{pr.1} is ill-posed 
because this is the property of the Cauchy problem for elliptic systems in ${\mathbb R}^n$  
(see, for instance \cite{Hd} or \cite[Ch. 1, \S 2]{mih}). 
Of course, in this case the boundary data should be taken independent on $t$. 
The Uniqueness Theorem clarify why the problem is ill-posed. The reason is the redundant  
data. Indeed, if $\Gamma$ has at least one interior point (on $\partial\Omega$), 
then taking a smaller set $\Gamma' \subset \Gamma$ we again obtain a problem 
with no more than one solution.

Another problem involves the initial data. 

\begin{problem} \label{pr.2}
Find 
a vector-function  
$u(x,\ t)\in 
[C^{2, 1}(\Omega_T) \cap C^{1,0}(\Omega_T \cup \overline{\Gamma_T}) \cap C _b 
(\overline{\Omega_T}) ]^n 
$ 
satisfying in $\Omega_T$ Lam\'e type equation {\rm (\ref{1.2})}, boundary conditions {\rm 
(\ref{2.1})}, {\rm (\ref{2.2})} and initial condition 
\begin{equation} \label{2.3}
u (x,0)=u^{(0)}(x),\ \ x\in \overline \Omega. 
\end{equation}
\end{problem}


Of course one should also take care on the compatibility of the data 
$u^{(0)}$, $u^{(1)}$, $u^{(2)}$: at least
\begin{equation} \label{eq.comp1}
u ^{(0)} (x) = u^{(1)} (x,0) \mbox{ on }  \overline \Gamma, 
\end{equation} 
and, if $u^{(0)}\in C^1 (\overline \Omega)$, even 
\begin{equation} \label{eq.comp2}
\sigma u^{(0)} (x) = u^{(2)} (x,0) \mbox{ on }  \overline \Gamma.
\end{equation}

%
%

The motivation of Problems \ref{pr.1} and \ref{pr.2} is transparent. 
The space $C_b (\Omega_T)$ is chosen because $u$ represents the ``velocity''. 
The first problem  describes the situation where for some reasons at each time $t\geq 0$ 
only part $\overline \Gamma$ of the solid surface $\partial  \Omega$ bounding the fluid is  
available for measurements. The second one describes the situation where the continuity up to  
$\partial \Omega_T$ is postulated, the ``velocity'' $u$ is known at every point  $x \in \overline 
\Omega$ at the initial time   $t=0$  but ґthe data on $\partial \Omega \setminus \overline 
\Gamma$ were lost for $t>0$. 

Clearly, Problem \ref{pr.2} has no more than one solution, too, if 
$\Gamma$ has at least one interior point (on $\partial\Omega$). 

We note that in classical theory of (initial and) boundary problems for the parabolic 
equation (\ref{1.2}), initial condition (\ref{2.3}) and boundary condition 
$\alpha u + \beta \sigma = u^{(3)}$ on the whole lateral surface ${\partial 
\Omega}_T$ of the cylinder $\Omega_T$ are usually considered. As a rule, such a problem is 
well-posed in proper spaces (H\"older spaces, Sobolev spaces etc.), see, for instance, 
\cite{lad}.

Let us show that Problem \ref{pr.2} is ill-posed, too.

\begin{example} \label{ex.1}
Let the Lam\'e coefficients $\mu$, $\lambda$ be constant and $a_j=0$, $0\leq j \leq n$. 

Take a cube $Q_n= \{0< x_j <1, \, 1\leq j \leq n\} $ as base $\Omega$ of the cylinder 
 $\Omega_T$. Let $\Gamma$ be the face $\{x_n = 0\}$ of the cube $Q_n$. Then 
 $\Gamma_T = Q_{n-1} \times (0,T)$ and the stress tensor $\sigma$ is given by 
the diagonal matrix with the non-zero entries
$$
\sigma_{j,j} =  \mu \frac{\partial}{\partial x_n} , 1\leq j \leq n-1, \, 
\sigma_{n,n} =  (2\mu + \lambda) \frac{\partial}{\partial x_n}. 
$$ 

Fix $N \in \mathbb N$ and consider  the sequence of functions $u (x,t,k,r)   \in [C^{\infty} 
({\mathbb R}^{n+1})]^n$ with the  components: 
 \begin{equation*}
 u_1 (x,t,k,r)=0, \, \dots , \, u_{n-1} (x,t,k,r)=0, \, u_n (x, t, 
 k,r)=\frac{e^{k^2(2\mu+\lambda)(t-r)+k x_n}}{k^N},
\end{equation*} 
depending on a parameter $0<r<+\infty$. Consider also the data $f (x,t,k,r)$,  $u^{(0)} 
(x,t,k,r)$, $u^{(1)} (x,t,k,r)$, $u^{(2)} (x,t,k,r)$ having the following components: 
\begin{equation*}
f_j (x,t,k,r)=0,  \, 1\leq j \leq n,  
\end{equation*}
\begin{equation*}
u^{(0)}_j (x,k,r)=0, \,   1\leq j \leq n-1,  \, 
u^{(0)}_n (x, k,r)=\frac{e^{-k^2 (2\mu+\lambda) r+k x_n}}{k^N}, 
\end{equation*}  
\begin{equation*}
u^{(1)}_j (x_1, \dots, x_{n-1},t,k,r)=0, \, 1\leq j \leq n-1,  
\end{equation*}  
\begin{equation*}
u^{(1)}_n (x_1, \dots, x_{n-1}, t, k,r)=\frac{e^{k^2 (2\mu+\lambda)(t-r)}}{k^N}, 
\end{equation*} 
\begin{equation*}
u^{(2)}_j (x_1, \dots, x_{n-1}, t, k,r)= 0, \, 1\leq j \leq n-1, 
\end{equation*} 
\begin{equation*}
 u^{(2)}_{n}  (x_1, \dots, x_{n-1}, 
t,k)= (2\mu+\lambda) \frac{e^{k^2(2\mu+\lambda)(t-T)}}{k^{N-1}}.
\end{equation*} 
Then, for $0<T<+\infty$, each function $u (x,t,k,T)$ is a solution 
to problem (\ref{1.2}), (\ref{2.1}), (\ref{2.2}), (\ref{2.3})  with 
the data $f (x,t,k,T)$,  $u^{(0)} (x,t,k,T)$, $u^{(1)} (x,t,k,T)$, $u^{(2)} (x,t,k,T)$.

It is clear, that compatibility conditions (\ref{eq.comp1}), (\ref{eq.comp2}) hold and
\begin{equation*}
f (x,t,k,T) \underset{k\rightarrow\infty}{\longrightarrow} 0 \mbox{ in } [C^\infty 
(\overline{\Omega_T})]^n, \quad  
u^{(0)} (x,k,T) \underset{k\rightarrow\infty}{\longrightarrow} 0 \mbox{ in } [C^\infty 
(\overline{\Omega})]^n, 
\end{equation*}
\begin{equation*}
u^{(1)} (x,t,k,T) \underset{k\rightarrow\infty}{\longrightarrow} 0 \mbox{ in } [C^s 
(\overline{\Gamma_T})]^n, \qquad 
u^{(2)} (x,t,k,T) \underset{k\rightarrow\infty}{\longrightarrow} 0 \mbox{ in } [C^{s} 
(\overline{\Gamma_T})]^n,  
\end{equation*}
if $N>2s+1$. On the other hand, for all $x_n>0$ and all  $N\in \mathbb N$ we have: 
\begin{equation*}
u_n (x, T, k)=\frac{e^{k^2 (2\mu+\lambda)(T-T)+k x_n}}{k^N}=\frac{e^{k 
x_n}}{k^N}\underset{k\rightarrow\infty}{\longrightarrow}+\infty.
\end{equation*}
Thus, there is no continuity with respect to the data and hence Problem 
\ref{pr.2} is ill-posed for $0<T<+\infty$. 
%

Let now $T=+\infty$. 
Then we may consider the data with a fixed $0<T_0<+\infty$:
$$
f (x,t,k,\infty)= 0 \in [C_b ( \overline{\Omega_T}]^n , \, 
u^{(0)} (x,k,\infty)= u^{(0)} (x,k,T_0)\in [C ( \overline \Omega)]^n, 
$$ 
$$
u_j^{(i)} (x,t,k,\infty)= 0, \, 1\leq j \leq n-1, \, 1\leq i \leq 2, 
$$
$$
u_n^{(1)} (x,t,k,\infty)= \left\{ 
\begin{array}{lll}
u_n^{(1)} (x,t,k, T_0), & t \leq T_0, \\
 \frac{1}{k^{N}} & t > T_0, \\
\end{array}
\right.
$$
$$
u_n^{(2)}(x,t,k,\infty) = \left\{ 
\begin{array}{lll}
u_n^{(1)} (x,t,k, T_0), & t \leq T_0, \\
 \frac{(2\mu+\lambda)}{k^{N-1}} & t > T_0. \\
\end{array}
\right.
$$
Obviously,  for $N\geq 2$,
\begin{equation*}
f (x,t,k,\infty) \underset{k\rightarrow\infty}{\longrightarrow} 0 \mbox{ in } 
[C_b (\overline{\Omega_T})]^n, \quad  
u^{(0)}, (x,k,\infty ) \underset{k\rightarrow\infty}{\longrightarrow} 0 \mbox{ in } [C^\infty 
(\overline{\Omega})]^n, 
\end{equation*}
\begin{equation*}
u^{(1)} (x,t,k,\infty) \underset{k\rightarrow\infty}{\longrightarrow} 0 \mbox{ in } [C^{1,0} 
(\overline \Gamma_T)\cap C _b 
(\overline{\Gamma_T})]^n, \,  
u^{(2)} (x,t,k,\infty) \underset{k\rightarrow\infty}{\longrightarrow} 0 \mbox{ in } [C _b 
(\overline {\Gamma_T})]^n. 
\end{equation*}
The Uniqueness Theorem \ref{t.1} for Problem \ref{t.1} implies that 
\begin{equation*}
u_n (x, t, k,\infty)= u_n (x, t, k,T_0) \mbox{ for } 0<t\leq T_0.
\end{equation*}
Then, for all $x_n>0$ and all  $2 \leq N\in \mathbb N$, we have 
$
\lim_{k\to +\infty}u_n (x, T_0, k,\infty) =+\infty.
$
Thus, if the data $f$, $u^{(0)}$, $u^{(2)}$, $u^{(2)}$ admits the solution to 
(\ref{1.2}) in $\Omega_{T}$ with boundary conditions  (\ref{2.1}), (\ref{2.2}) 
for $T=+\infty$ and the initial condition (\ref{2.3})
 then there is no continuity with respect to the data in 
the chosen space. Otherwise there is no solutions to the problem 
for some data in the data's spaces. In any case, the problem is ill-posed for $T=+\infty$, too.
\end{example}

As both Problems \ref{pr.1} and  \ref{pr.2} are ill-posed, we will not study Problem 
\ref{pr.2} because in addition to (\ref{1.2})-(\ref{2.2}) to investigate it 
one needs to know also the data related to initial condition (\ref{2.3}). 
Besides, we will consider the case $0<T<+\infty$ only.


\section{Solvability Conditions}

From now on we will study Problem \ref{pr.1} under the assumption that 
its data belong to H\"older spaces (cf., \cite[Ch. 1, \S 1]{frid} for other 
boundary problems for parabolic equations). We recall that  
a function $u(x)$, defined on a set $M\in{\mathbb R}^m$, is called \emph{H\"older  continuous 
with a power $0<\lambda<1$} on $M$,  if there is such a constant  $C>0$ that 
\begin{equation} \label{1.1}
|u(x)-u(y)|\leq C |x-y|^\lambda \mbox{ for all } x, y \in M 
\end{equation}
($|x-y| = \sqrt{\sum_{j=1}^m (x_j-y_j)^2}$ being Euclidean distance between points 
$x$ and $y$ in ${\mathbb R}^m$). Let $C^{\lambda} (\overline{\Omega_T})$ stand for the set of 
H\"older  continuous functions with a power $\lambda$ over $\overline{\Omega_T}$. Besides, let 
$C^{1+\lambda, \lambda} (\overline{\Omega_T})$ be the set of H\"older  continuous functions with 
a power $\lambda$ over $\overline{\Omega_T}$, having H\"older  continuous derivatives 
$u_{x_i}$, $1\leq i \leq n$, with the same power in $\overline{\Omega_T}$.

We choose a set $\Omega^+$ in such a way that the set $D=\Omega\cup\Gamma\cup\Omega^+$ would 
be a bounded domain with piece-wise smooth boundary. It is possible since  $\Gamma$ is an open 
connected set. It is convenient to set $\Omega^- = \Omega$. For a function $v$ on $D_T$ we denote 
by  $v^+$ its restriction to $\Omega^+$ and, similarly, we denote by  $v^-$ its restriction to 
$\Omega$. It is natural to denote limit values of  $v^\pm$ on $\Gamma_T$, when they are defined, 
by  $v^\pm_{|\Gamma_T}$. 

\begin{theorem}[Solvability criterion] \label{t.2}
Let  $\Gamma \in C^{1+\lambda}$, 
$$
f\in [C^{\lambda}(\overline{\Omega_T})]^n, \, u_1\in [C^{1+\lambda,\lambda}
(\overline{\Gamma_T})]^n, \, u_2\in [C^{\lambda}(\overline{\Gamma_T})]^n.
$$ 
Problem \ref{pr.1} 
is solvable 
if and only if there is a vector-function $F\in [C^{2,1}(D_T)]^n$ satisfying the following 
conditions:
\begin{enumerate}
\item[{\rm 1)}] $L_{n+1} F=0$ in $D_T$,
\item[{\rm 2)}] $F=G_{\Omega, 0} (f)+V_{\overline{\Gamma},0} (u_2) +
W_{\overline{\Gamma},0} (u_1)$ in $\Omega^+_T$.
\end{enumerate}
\end{theorem}

\textbf{Proof.}\emph{ \textrm{Necessity.}} Let a function  $u(x, t)\in [C^{2,1}(\Omega_T)\cap 
C^{1,0}(\Omega_T\cup\overline{\Gamma_T})]^n$ satisfies (\ref{1.2}), (\ref{2.1}), (\ref{2.1}).
Consider the function 
\begin{equation*}
F=G_{\Omega, 0} (f)+V_{\overline{\Gamma},0} (u_2) + 
W_{\overline{\Gamma},0} (u_1) -\chi_{\Omega_T} u.
\end{equation*} 
in the domain $D_T$, where $\chi_{M}$ is a characteristic function of the set 
$M \subset {\mathbb R}^{n+1}$. By the very construction condition 2) is fulfilled for it. 

Clearly, the function  $u(x,t)$ belongs to the space $[C^{1,2}(\overline{\Omega'_T})]^n$ for each 
cylindrical domain $\Omega'_T$ with such a base $\Omega'$ that $\Omega' \subset \Omega$ and 
$\overline{\Omega'}\cap\partial \Omega \subset \Gamma$. Besides,  
$L_{n+1}u =f \in [C^{\lambda}(\overline{\Omega_T'})]^n$. Without loss of the generality we may assume 
that the interior part $\Gamma' $ of the set  $\overline{\Omega'}\cap\partial \Omega$ 
is non-empty.

We note that  $\chi_{\Omega_T} u=\chi_{\Omega'_T} u$ in $D'_T$, where $D'= \Omega' \cup 
\Gamma' \cup \Omega^+$. Then using Lemma \ref{l.1} we obtain:
\begin{equation} \label{3.5}
F=G_{\Omega \setminus \overline{\Omega'}, 0} (f)+V_{\overline{\Gamma}\setminus \Gamma',0} (u_2) + 
W_{\overline{\Gamma}\setminus \Gamma',0} (u_1) - I_{\Omega',0} (u) \mbox{ in } D'_T.
\end{equation} 

Arguing as in the proof of Theorem  \ref{t.1} we conclude that  
each of the integrals in the right hand side of (\ref{3.5}) satisfies 
homogeneous Lam\'e type equation outside the corresponding integration set. 
In particular, we see that $L_{n+1}F =0$ in $D'_T$. Obviously, for any point $(x, t)\in D_T$ 
there is a domain $D'_T$ containing $(x,\ t)$. That is why  $L_{n+1}F =0$ in $D_T$, and 
hence $F$ belongs to the space $[C^{2,1}(D_T)]^n$. Thus this function satisfies 
condition 1), too. 

\emph{\textrm{Sufficiency.}} Let there be a function $F \in  [C^{2,1}(D_T) ]^n$, satisfying 
conditions 1) and 2) of the theorem. Consider on the set $D_T$ the function 
\begin{equation} \label{eq.sol}
U=G_{\Omega, 0} (f)+V_{\overline{\Gamma},0} (u_2) +
W_{\overline{\Gamma},0} (u_1) - F .
\end{equation} 
As  $f \in [C^\lambda(\overline{\Omega_T})]^n$ then the results of \cite[Ch. 1, \S 3]{frid} imply 
\begin{equation} \label{eq.G1} 
G_{\Omega, 0} (f) \in [C^{2,1}(\Omega^\pm_T) \cap 
C^{1,0}(D_T)  \cap C (\overline D_T)]^n
\end{equation} 
and, moreover, 
\begin{equation} \label{eq.G2}  
L_{n+1} G^-_{\Omega, 0} (f) = f \mbox{ in } \Omega_T, \quad L_{n+1} G^+_{\Omega, 0} (f) = 0 
\mbox{ in }  \Omega^+_T.
\end{equation} 
Since  $u_2 \in [C^{\lambda} (\overline{\Gamma_T})]^n$ then the results of  \cite[Ch. 5, \S 
2]{frid}yield 
\begin{equation} \label{eq.V1}  
V_{\overline \Gamma, 0} (u_2) \in [C^{2,1}(\Omega^\pm_T) \cap 
C^{1,0}((\Omega^\pm\cup \Gamma)_T) 
\cap C(\overline{D_T}   \setminus (\partial \Gamma)_T )]^n, 
\end{equation} 
\begin{equation}  \label{eq.V2}   
L_{n+1} V_{\overline \Gamma, 0} (u_2) = 0 \mbox{ in } \Omega_T \cup \Omega^+_T .
\end{equation} 
On the other hand, the behavior of the Double Layer Potential $W_{\overline \Gamma, 0} (u_1)$ 
is similar to the behavior of the normal derivative of Single 
Layer Potential $V_{\overline \Gamma, 0} (u_1)$. Hence 
\begin{equation} \label{eq.W1}  
W_{\overline \Gamma, 0} (u_1) \in [C^{2,1}(\Omega^\pm_T) \cap 
C(\overline{\Omega^\pm_T} \setminus (\partial \Omega^\pm \setminus \Gamma)_T)]^n,
\end{equation}
\begin{equation} \label{eq.W} 
L_{n+1} W_{\overline \Gamma, 0} (u_1) = 0 \mbox{ in } \Omega_T \cup \Omega^+_T.
\end{equation}   

\begin{lemma} \label{l.1}
Let $S \subset \overline \Gamma \in C^{1+\lambda}$.  
If $u_1\in [C^{1+\lambda,\lambda}(\overline{\Gamma_T})]^n$, then the potential $W^-_{\overline 
\Gamma, 0} (u_1)$ belongs to the space $[C^{1, 0}(\Omega_T\cup S_T)]^n$ if and only if 
$W^+_{\overline \Gamma, 0} (u_2) \in [C^{1,0}(\Omega^+_T\cup S_T)]^n$.
\end{lemma}

\textbf{Proof.} It is similar to the proof of 
the analogous lemma for Newton Double Layer Potential (see, for instance, \cite[lemma 1.1]{Sh1}).
Actually, one needs to use Lemma  \ref{l.1} instead of the standard Green formula for the 
Laplace operator.
\hfill $\square$

Since $F \in [C^{1,0} (D_T)]^n$ then it follows from the discussion above that 
$W^+_{\overline \Gamma, 0} (u_2) \in [C^{1,0}((\Omega^+\cup \Gamma)_T)]^n$.
Thus, formulas (\ref{eq.sol})--(\ref{eq.W}) and Lemma \ref{l.1} imply that  
\begin{equation*} 
U \in [C^{2,1}(\Omega^\pm_T ) 
\cap C^{1,0}((\Omega^\pm \cup \Gamma)_T) 
\cap C(\overline{\Omega^\pm_T} \setminus (\partial \Omega \setminus \Gamma)_T)]^n ,
\end{equation*} 
\begin{equation*} 
L_{n+1} U = \chi_{D_T} f \mbox{ in } \Omega_T \cup \Omega^+_T.
\end{equation*}   
In particular, (\ref{1.2}) is fulfilled for $U^-$.

Let us show that the function $U^-$ satisfies (\ref{2.1}) and (\ref{2.2}). 

Since  $F \in  [C^{1,0} (D_T)]^n $ we see that $\partial ^\alpha F^- = \partial ^\alpha F^+$ on 
$\Gamma_T$ for $\alpha \in {\mathbb Z}_+$ with $|\alpha|\leq 1$ and 
\begin{equation*} 
\partial ^\alpha F^+_{|\Gamma_T} = \left( \partial ^\alpha  G^+_{\Omega, 0} (f)+ \partial ^\alpha  
V^+_{\overline{\Gamma},0} (u_2) + \partial ^\alpha W^+_{\overline{\Gamma},0} 
(u_1)\right)_{|\Gamma_T}.
\end{equation*}  
It follows from formulas (\ref{eq.G1}) and (\ref{eq.V1}) that the  Parabolic Volume Potential and 
the Single Layer Parabolic Potential are continuous if the point  $(x,t)$ passes over the surface 
$\Gamma_T$. Then 
\begin{equation*} 
U^-_{|\Gamma_T}=
W^-_{\overline{\Gamma},0} (u_1)_{|\Gamma_T} - W^+_{\overline{\Gamma},0} (u_1) _{|\Gamma_T} = u_1.
\end{equation*}  
because of the theorem on jump behavior of the Parabolic Double Layer Potential 
(see, for instance, \cite[Ch. 5, \S 2, theorem 1]{frid}), i.e. equality (\ref{2.1}) is valid 
for $U^-$ .  

Formula (\ref{eq.G1}) means that that the surface stress of the Parabolic Volume Potential 
is continuous if the point  $(x,t)$ passes over the surface $\Gamma_T$. Therefore  
\begin{equation}  \label{eq.U1}
(\sigma U)^-_{|\Gamma_T}=\left( \sigma  
V^-_{\overline{\Gamma},0} u_2 \right) _{|\Gamma_T} - \left(  \sigma 
 V^+_{\overline{\Gamma},0} u_2\right)_{|\Gamma_T} + 
\left( \sigma W^-_{\overline{\Gamma},0} u_1 \right)_{|\Gamma_T} - \left(
\sigma  W^+_{\overline{\Gamma},0} u_1 \right)_{|\Gamma_T} .
\end{equation}  
By theorem on jump behavior of the stress of the Parabolic Single Layer Potential
(see, for instance, \cite[Ch. 3, \S 10, theorem 10.1]{landis}) 
\begin{equation}  \label{eq.U2}
\left(\sigma V^-_{\overline{\Gamma},0} u_2 \right) _{|\Gamma_T} - \left( \sigma 
V^+_{\overline{\Gamma},0} u_2 \right)_{|\Gamma_T} = u_2 .
\end{equation}  

Finally, we need the following lemma which is an analogue of the famous Theorem  
on jump behavior of the normal derivative of the Newton's Double Layer Potential.

\begin{lemma} \label{l.2} Let $\Gamma \in C^{1+\lambda}$ and $u_2 \in [C ^\lambda (\overline{\Gamma_T})]^n $.
If $W^-_{\overline{\Gamma},0} (u_1) \in [C^{1, 0}((\Omega \cup \Gamma)_T)]^n $ 
or $W^+_{\overline{\Gamma},0} (u_1) \in [C^{1, 0}((\Omega^+ \cup 
\Gamma)_T)]^n $ then   
\begin{equation} \label{eq.Wnu} 
\left( \sigma W^-_{\overline{\Gamma},0} u_1 - 
\sigma W^+_{\overline{\Gamma},0} u_1\right)_{|\Gamma_T} =0. 
\end{equation}  
\end{lemma}

{\bf Proof.} Really, let, for instance,  $W^-_{\overline{\Gamma},0} (u_1) \in [C^{1, 
0}((\Omega \cup \Gamma)_T)]^n $.  Then using  Lemma \ref{l.1} we obtain 
$W^+_{\overline{\Gamma},0} u_1 \in [C^{1, 0}((\Omega^+ \cup \Gamma)_T)]^n $ and   
$\left(\sigma W^\pm_{\overline{\Gamma},0} (u_1)\right)_{|\Gamma_T} \in 
[C(\Gamma_T)]^n$.

Let $\phi \in [C^\infty_0 (D_T)]^n $ be a function with compact support in $D_T$. Then 
formulas (\ref{eq.Green.0})--(\ref{eq.Green.d})
yield:
\begin{equation} \label{eq.Wnu1}
\int_{\Gamma_T} \phi^* \left( \sigma W^-_{\overline{\Gamma},0} u_1 - 
\sigma W^+_{\overline{\Gamma},0} u_1\right) ds (x) dt = 
\end{equation}
\begin{equation*}
\int_{\Omega_T\cup\Omega^+_T} \phi^* ( {\mathcal L}_n +a) W_{\overline{\Gamma},0} u_1 dx dt  
 + \int_{T_1}^{T_2} {\mathfrak D}_{\Omega\cup\Omega^+} (W_{\overline{\Gamma},0} u_1,\phi)dt =
 \end{equation*}
\begin{equation*}
\int_{\Omega_T\cup\Omega^+_T} \phi^* \Big( \frac{\partial}{\partial t}  -A + a\Big)
W_{\overline{\Gamma},0} u_1 dx dt
+ \int_{T_1}^{T_2} {\mathfrak D}_{\Omega\cup\Omega^+} (W_{\overline{\Gamma},0} u_1,\phi)dt
 \end{equation*}
because $L_{n+1} W^\pm_{\overline{\Gamma},0} u_1 =0$ in $\Omega^\pm$ according to  (\ref{eq.W}).

Again, integrating by parts and using formulas (\ref{eq.Green.0})--(\ref{eq.Green.d}) and 
Theorem on jump behavior of the Parabolic Double Layer 
Potential, 
we see that 
\begin{equation} \label{eq.Wnu2}
\int_{\Omega_T\cup\Omega^+_T} \phi^* \Big( \frac{\partial}{\partial t} -A +a\Big)  
W_{\overline{\Gamma},0} u_1 dx dt   +  \int_{T_1}^{T_2} {\mathfrak D}_{\Omega\cup\Omega^+} (W_{\overline{\Gamma},0} u_1,\phi)dt = 
  \end{equation}
\begin{equation*} 
-\int\limits_{\Omega_T\cup\Omega^+_T} \Big(\frac{\partial \phi}{\partial t} \Big)^* 
 W_{\overline{\Gamma},0} u_1 dx dt   - \int\limits_{\Omega_T\cup\Omega^+_T} (({\mathcal L}_n +A^*)
 \phi)^*  W_{\overline{\Gamma},0} u_1 dx dt + 
 $$
 $$
 \int_{\Gamma_T}  (\tilde \sigma \phi)^*  (W^-_{\overline{\Gamma},0} u_1  - 
 W^+_{\overline{\Gamma},0} u_1 )ds (x) dt =
 \end{equation*}
\begin{equation*} 
 \int_{\Gamma_T}  (\tilde \sigma \phi)^*  u_1 ds (x) dt 
- \int_{\Omega_T\cup\Omega^+_T} (L_{n+1}^*\phi)^* W_{\overline{\Gamma},0} u_1  dx dt  .
 \end{equation*}
But the kernel $\Phi (x,y, t ,\tau)$ is a fundamental solution of the backward 
parabolic operator $L_{n+1}^*$ with respect to variables $(y,\tau)$. Hence  
\begin{equation*}
\int_{D_T} (L_{n+1}^*\phi (x,t))^* \Phi (x,y, t ,\tau) dx d t = \phi^* (y,\tau) , 
\quad (y,\tau) \in D_T.
 \end{equation*}
Then the type of the singularity of the fundamental solution allows us to apply Fubini Theorem 
and to conclude that 
\begin{equation} \label{eq.Wnu3}
\int_{\Omega_T\cup\Omega^+_T} (L_{n+1}^*\phi)^* W_{\overline{\Gamma},0} u_1 dx dt = 
   \end{equation}
\begin{equation*}   
 \int_{\Gamma_T} \tilde \sigma \int_{D_T} (L_{n+1}^*\phi (x,t))^* 
 \Phi (x,y, t ,\tau) dx d t  \, u_1  \, ds (y) d\tau  =  \int_{\Gamma_T}  (\tilde 
 \sigma \phi)^*  u_1 ds (y) d\tau.
  \end{equation*}
Finally, formulas (\ref{eq.Wnu1})- (\ref{eq.Wnu3}) imply that   
\begin{equation*} 
\int_{\Gamma_T} \phi^* \left( \sigma W^-_{\overline{\Gamma},0} u_1  - 
\sigma  W^+_{\overline{\Gamma},0} u_1\right) ds = 0
 \end{equation*}
for all $\phi \in [C^\infty_0 (D_T)]^n$. As such functions are dense in the Lebesgue space  
$[L^1(K)]^n$ for any compact $K \subset \Gamma_T$ then formula (\ref{eq.Wnu}) holds true. 
 \hfill $\square$

Now using lemma \ref{l.2} and formulas (\ref{eq.U1}), (\ref{eq.U2}), we conclude that  
$(\sigma U)^-_{|\Gamma_T}= u_2$,
i.e. (\ref{2.2}) is fulfilled for $U^-$.

Thus, function $ u (x,t)=U^-(x, t)$ satisfies conditions (\ref{1.2})--(\ref{2.2}). 
The proof is complete. \hfill $\square$

It follows from formula (\ref{eq.sol}) that properties of a solution 
to Problem \ref{pr.1} depend on  properties of the extension $F$ of 
the sum of the parabolic potentials, described in Theorem \ref{t.2}.

\begin{corollary} \label{c.2}
Let $S \subset \partial \Omega \setminus \Gamma$.  
Under the hypotheses of Theorem \ref{t.2}, Problem
\ref{pr.1} is solvable in the space
\begin{equation*}
[C^{2,1}(\Omega_T) \cap C^{1,0}(\Omega_T \cup \overline{\Gamma_T}) \cap 
C(\overline{\Omega_T}\setminus S_T)]^n
\end{equation*}
if and only if there exists a function 
\begin{equation*}
F\in [C^{2,1}(D_T) \cap C^{1,0}(\Omega_T \cup 
\overline{\Gamma_T}) \cap C(\overline{\Omega_T}\setminus S_T)]^n, 
\end{equation*}
satisfying conditions 1) and 2) of Theorem \ref{t.2}.
\end{corollary}

In particular, if $S = \emptyset$ then corollary \ref{c.2} gives criterion for 
the existence of solution to Problem \ref{pr.1} in the space $[C (\overline{\Omega_T})]^n$.

We note that Theorem \ref{t.2} is an analogue of Theorem by Aizenberg and Kytmanov 
\cite{AKy}) describing solvability conditions of the Cauchy problem for the Cauchy--Riemann 
system (cf. also \cite{Sh1} in the Cauchy Problem for Laplace Equation or  
 \cite{Tark36} in the  Cauchy problem for general elliptic systems).
Formula (\ref{eq.sol}), obtained in the proof of Theorem  \ref{t.2}, gives 
the unique solution to Problem  \ref{pr.1}. Clearly, if we will be able to write  
the extension $F$ of the sum of potentials $G_{\Omega, 0} (f)+V_{\overline{\Gamma},0} (u_2) +
W_{\overline{\Gamma},0} (u_1)$ from $\Omega^+$ onto $D_T$ as a series with respect to 
special functions or a limit of parameter depending integrals then we will get  
a Carleman type formula for 
solutions to Problem \ref{pr.1} (cf. \cite{AKy}). However this is a topic 
for another paper. In the sequel we will discuss polynomial and 
formal solutions for operators with constant coefficients only. 

\section{Polynomial solutions and dense solvability}

It is not difficult to prove dense solvability of Problem \ref{pr.1} 
in the case where $\Gamma$ is an open connected set of the hyperplane $\{x_n =0\}$.

\begin{lemma} \label{l.poly} If $\Gamma$ is an open connected set if the hyperplane $\{x_n =0\}$ the 
 Problem \ref{pr.1}  is densely solvable.
\end{lemma}

{\bf Proof.}  First let us prove that if in this case  
the data of Problem \ref{pr.1} are polynomials then the problem is solvable and its solution 
is a polynomial. 

Indeed, Problem \ref{pr.1} is easily can be reduced to the following one (see Example 
\ref{ex.1}):  
\begin{equation} \label{1.2a}
L_{n+1}v = g \mbox{ in } \Omega_T 
\end{equation}
\begin{equation} \label{2.1a}
 v(x_1, \dots, x_{n-1}, 0,t)=0 \mbox{ on } \overline{\Gamma_T}, 
\end{equation}
\begin{equation} \label{2.2a}
\mu \frac{\partial v_j}{\partial x_n} (x_1, \dots, x_{n-1}, 0,t) =0 \mbox{ on } 
\overline{\Gamma_T}, 1\leq j \leq n-1, 
\end{equation}
\begin{equation} \label{2.2a.n}
(2\mu +\lambda) \frac{\partial v_n}{\partial x_n} (x_1, \dots, x_{n-1}, 0,t) =0 \mbox{ on } 
\overline{\Gamma_T}, 1\leq j \leq n-1, 
\end{equation}
with 
\begin{equation*}
g (x,t)= f (x,t)- (L_{n}  u_1) (x_1, \dots, x_{n-1},t) 
- x_n J (\mu,\lambda) (L_{n} u_2) ( x_1, \dots, x_{n-1},t) .
\end{equation*} 
where 
$J (\mu,\lambda)$ is the diagonal matrix with the components 
$$
J_{j,j} (\mu,\lambda) = \mu^{-1}, \, 1\leq j \leq n-1, \, J_{n,n} (\mu,\lambda) = 
(2\mu + \lambda)^{-1}.
$$
Besides, 
$
u (x,t)= v (x,t)+ u_1 ( x_1, \dots, x_{n-1},t) + J (\mu,\lambda)
 x_n u_2 (x_1, \dots, x_{n-1},t) .
$  

Now consider data $g^{(j,\alpha)} (x,t) = t^j x^\alpha$ with a multi-index $\alpha \in {\mathbb Z}_+^n$. 

If $0 \leq \alpha_1 + \dots \alpha_{n-1} \leq 1$, we easily obtain (unique) polynomial solutions
\begin{equation} \label{eq.w}
v^{(j,\alpha)} (x,t) = x_1 ^{\alpha_1} \cdots x_{n-1} ^{\alpha_{n-1}} 
w^{(j,\alpha_n)} (x_n,t), \, \alpha_n, j \in {\mathbb Z}_+,  
\end{equation} 
to  problem (\ref{1.2a})--(\ref{2.2a}) where 
\begin{equation*}
w^{(0,k)} (y,t)= - \frac{y^{k+2 } k!}{(k+2)!}, 
\quad w^{(1,k)} (y,t) = - \frac{t y^{k+2 } k!}{(k+2)!} - 
\frac{ y^{k+4} k!}{(k+4)!}, \, k \in {\mathbb Z}_+, \, y \in {\mathbb R}
\end{equation*}
and, by the induction with respect to 
$j \in {\mathbb Z}_+$, 
\begin{equation} \label{eq.ww}
w^{(j,k)} (y,t) = - \sum_{\mu=0}^j \frac{t^{j-\mu} y^{k+2\mu+2} k! j ! }
{(k+2\mu + 2)! (j-\mu)!}, \, k \in {\mathbb Z}_+, \, y \in {\mathbb R}.
\end{equation}

To finish the arguments we use the induction with respect to $|\alpha '| \in {\mathbb Z}_+$ where  
$\alpha' = (\alpha_1 , \dots, \alpha_{n-1}) \in {\mathbb Z}_+^{n-1}$. Namely, let 
for $s \geq 2$ and all $\alpha'$ with $|\alpha '| =s$ the solutions to the problem are polynomial. 
If  $|\alpha '| =s+1$ then 
\begin{equation*}
L_{n+1}  \left( x_1 ^{\alpha_1} \cdots x_{n-1} ^{\alpha_{n-1}} 
w^{(j,\alpha_n)} (x_n,t)\right)  =  t^j x^\alpha - w^{(j,\alpha_n)} (x_n,t)  \Delta_{n-1} \left( x_1 
^{\alpha_1} \cdots x_{n-1} ^{\alpha_{n-1}} \right).
\end{equation*}
Clearly, the degree of the polynomial $p_{j,\alpha} 
(x,t)= w^{(j,\alpha_n)} (x_n,t)  \Delta_{n-1} \left( x_1 
^{\alpha_1} \cdots x_{n-1} ^{\alpha_{n-1}} \right)$ with respect to $x' \in {\mathbb R}^{n-1}$ 
equals to $s-1$. Then, by the induction, problem (\ref{1.2})--(\ref{2.2}) with data 
$p_{j,\alpha} (x,t)$ admits a polynomial solution, say, $r_{j,\alpha}  (x,t)$. Therefore the solution $v^{(j,\alpha)} (x,t)$ to problem (\ref{1.2})--(\ref{2.2}) with data $g^{(j,\alpha)} (x,t)= 
t^j x^\alpha$, $|\alpha '| =s+1$,  is given as follows: 
\begin{equation*}
v^{(j,\alpha)} (x,t) = x_1 ^{\alpha_1} \cdots x_{n-1} ^{\alpha_{n-1}} 
w^{(j,\alpha_n)} (x_n,t)  + r_{j,\alpha} (x,t),
\end{equation*} 
i.e. it is a polynomial, too. 

Now Problem  \ref{pr.1} with zero boundary data in the 
 case $\Gamma \subset \{x_n =0\}$  is densely solvable because any continuous function  
$g$ on the compact set $\overline{\Omega_T}$ can be approximated by polynomials.
But the reducing to zero boundary data was organized in such a way that 
 one easily sees, in this case Problem  \ref{pr.1} is densely solvable for 
 non-zero boundary data, too. 
\hfill $\square$

The dense solvability of Problem \ref{pr.1} in 
 general setting 
is natural to expect if the set $\partial \Omega \setminus \overline \Gamma$ has at least 
one interior point in $\partial\Omega$ (cf. \cite{ShTaLMS} in the Cauchy 
Problem for elliptic equations). 

Finally, we note that polynomial solutions indicated in the proof of Lemma \ref{l.poly} can be 
used in order to construct formal solutions to Problem \ref{pr.1}. 

\section{Basis with double orthogonality}

Denote by $\{h_\nu ^{(s)}\}$ the set of harmonic homogeneous polynomials 
(spherical harmonics) forming an orthonormal basis in $L^2 (\partial B(0,1))$;
here $\nu$ is the 
degree of the homogeneity and $1\leq s \leq J(s)$ where  $J (s) = \dfrac{(2n + 2s - 2)(2n + s - 
3)!}{s! (2n-2)!}$ (see \cite{So}).

\begin{lemma} \label{l.bdo}
Let $\mu=1$, $\lambda=-1$. Then the polynomials
$$
H^{(s)}_{0,\nu} (x,t)= h^{(s)}_\nu (x), \nu \geq 0, 
$$
$$
 H^{(s,i)}_{N,0} (x,t) = \sum_{j=0}^N \frac{t^j 
x_i^{2N-2j} (2N)!}{j! (2N-2j)!} , \, 1\leq i \leq n , N \geq 1, 
$$
$$
 H^{(s)}_{N,\nu} (x,t) = \sum_{j=0}^N \frac{t^j |x|^{2N-2j} h_\nu^{ (s)} (x) (2N)!! 
 (n+2N-2+2\nu) !!}{j! (2N-2j)!! (n+2N-2j+2\nu)!!} , \nu \geq 1, N \geq 1
$$
are solutions to the heat equations in ${\mathbb R}^{n+1}$.

Besides, $H^{(s)}_{N,\nu}$ and $H^{(p)}_{M,\mu}$ are $L^{2} (B (0,R)_T)$-orthogonal 
for all $R>0$ and $T>0$ if $(\nu, s) \ne (\mu,p)$.
\end{lemma}

{\bf Proof.} Indeed, for $\nu \geq 0$, $N\geq 1$, and $1\leq i \leq n$ we have:
$$
L_{n+1} H^{(s)}_{0,\nu} = -\Delta_n h^{(s)}_\nu , 
$$
$$
L_{n+1} H^{(s,i)}_{N,0}  =  \sum_{j=1}^N \frac{t^{j-1} 
x_i^{2N-2j} (2N)!}{(j-1)! (2N-2j)!} - \sum_{j=0}^{N-1} \frac{ t^{j} 
x_i^{2N-2j-2} (2N)! }{j! (2N-2j-2)!} =0.
$$

On the other hand, 
$$
\frac{\partial |x|^{2k}}{\partial x_j} = 2 k x_j |x|^{2k-2} , \, 
\frac{\partial^2 |x|^{2k}}{\partial x^2_j} = 2 k  |x|^{2k-2} + 2 k (2k-2) x^2_j |x|^{2k-4}
$$
$$
\Delta_n |x|^{2k} = 2k (n+2k-2) |x|^{2k-2}, k \geq 1. 
$$
Hence, for $\nu \geq 1$, $k \geq 1$, 
$$
\Delta_n (|x|^{2k} h_\nu (x)) = 2k (n+2k+2\nu-2) |x|^{2k-2} h_\nu ^{(s)} (x)
$$
because of Euler's formula 
$$
\sum_{j=1}^n x_j \frac{\partial h_\nu^{(s)} (x)}{\partial x_j} = \nu h_\nu ^{(s)} (x). 
$$

Consider the polynomial $H=\sum_{j=0}^N  c_j  t^j |x|^{2N-2j} h_\nu^{ (s)} (x) $
with constants $c_j$. Then 
$$
L_{n+1} H = \sum_{j=1}^N  c_j  j t^{j-1} |x|^{2N-2j} h_\nu^{ (s)} (x) 
- 
$$
$$
\sum_{j=0}^{N-1}  c_j   t^{j} 
(2N-2j) (n+2N-2j+2\nu-2) |x|^{2N-2j-2} h_\nu^{ (s)} (x) =
$$
$$
\sum_{j=1}^N  t^{j-1} |x|^{2N-2j} h_\nu^{ (s)} (x) [jc_j - c_{j-1} (2N-2j+2)(n+2N-2j+2\nu) ]
$$
Thus, we get a recurrent formula 
$$
c_j = j^{-1}c_{j-1} (2N-2j+2)(n+2N-2j+2\nu) 
$$
for the coefficients in the case $L_{n+1} H=0$. Choosing $c_0=1$ we easily obtain 
$$
c_j =  \frac{(2N)!! (n+2N-2+2\nu) !!}{j! (2N-2j)!! (n+2N-2j+2\nu)!!}.
$$

Finally the statement on $L^{2} (B (0,R)_T)$-orthogonality follows from Fubini 
Theorem and the homogeneity of the polynomials $h^{(s)}_\nu (x)$.
$\hfill \square$

This lemma suggests us to consider a function 
$$
H_\nu(x,t)=\phi(|x|,t) h_\nu ^{(s)} (x).
$$
Easily
$$
\frac{\partial\phi (|x|,t)}{\partial x_j} = \frac{x_j}{|x|}  \frac{\partial\phi }{\partial z} (|x|,t),
$$
$$
\frac{\partial^2\phi (|x|,t)}{\partial x^2_j} = \frac{1}{|x|}  \frac{\partial\phi }{\partial z} (|x|,t) - 
\frac{x^2_j}{|x|^3}  \frac{\partial\phi }{\partial z} (|x|,t) +  \frac{x^2_j}{|x|^2}  \frac{\partial^2 \phi }{\partial^2 z} (|x|,t),
$$
$$
\sum_{j=1}^n x_j \frac{\partial h_\nu }{\partial x_j} = \nu h_\nu.
$$
Hence
$$
L_{n+1} H_{\nu} (x,t) =\Big( 
\frac{\partial\phi}{\partial t}(|x|,t)  - \frac{\partial^2\phi}{\partial z^2}(|x|,t) - \frac{(n+2\nu-1)}{|x|}\frac{\partial\phi}{\partial 
z}(|x|,t) \Big) h_\nu ^{(s)} (x) =0
$$ 
if $\phi (z,t)$ is a solution to parabolic equation 
$$
\frac{\partial\phi}{\partial t}(z,t)  - \frac{\partial^2\phi}{\partial z^2}(z,t) - \frac{(n+2\nu-1)}{z}\frac{\partial\phi}{\partial 
z}(z,t) =0
$$
for $z>0$, $t>0$

\bigskip

\emph{The work was supported by RFBR grant 11-01-00852a.
}


\end{document}